\documentclass[10pt]{article}
\usepackage{amsfonts}
\usepackage{latexsym}
\usepackage{cite}
\usepackage{amsmath,amsfonts,latexsym,amssymb}
\usepackage[mathscr]{eucal}
\usepackage{cases,color}
\usepackage{amsthm}

\usepackage[bf,small]{caption2}
\usepackage{float}
\usepackage{graphicx}
\usepackage{amsmath}
\usepackage{amssymb}
\usepackage[all]{xy}

\newtheorem{theorem}{theorem}[section]

\newtheorem{lem}[theorem]{Lemma}

\newtheorem{thm}[theorem]{Theorem}

\begin{document}

\title{\textbf{Lifting homeomorphisms and finite abelian branched covers of the 2-sphere}}
\author{\Large Haimiao Chen}
\date{}
\maketitle

\begin{abstract}
  We completely determine finite abelian regular branched covers of the 2-sphere $S^2$ with the property that each homeomorphism of $S^2$ preserving the branching set can be lifted.

\end{abstract}

\section{Result and method}

As pointed out by Birman and Hilden \cite{BH17} (see also \cite{AMO21,GW17,MW21}), an interesting problem is to find branched covers of a surface with the property that each homeomorphism of the base surface preserving the branching set can be lifted.

Fix a set $\mathfrak{B}=\{\mathfrak{x}_1,\ldots,\mathfrak{x}_{n}\}\subset S^2$, with $n\ge 2$. Let $\Sigma_{0,n}=S^2\setminus\mathfrak{B}$. Let $x_i\in H_1(\Sigma_{0,n};\mathbb{Z})$ denote the homology class of a small loop enclosing $\mathfrak{x}_i$.
Given an abelian group $A$, a regular $A$-cover $\pi:\Sigma\to S^2$ with branching set $\mathfrak{B}$ is determined by an epimorphism $\phi_\pi: H_1(\Sigma_{0,n};\mathbb{Z})\twoheadrightarrow A$, which satisfies $\phi_\pi(x_1)+\cdots+\phi_\pi(x_n)=0$.
Call two such covers $\phi:\Sigma\to S^2$ and $\phi':\Sigma'\to S^2$ {\it equivalent} if there exist homeomorphisms $\tilde{g}:\Sigma\to\Sigma'$, $g:S^2\to S^2$ such that $g(\mathfrak{B})=\mathfrak{B}$ and $g\circ\pi=\pi'\circ\tilde{g}$.

We extend the results of \cite{AMO21,GW17} by showing
\begin{thm}\label{thm:main}
Suppose $A$ is a finite abelian $p$-group with exponent $p^k$, and $\pi:\Sigma\to S^2$ is a regular $A$-cover with branching set $\mathfrak{B}$ such that each homeomorphism of $S^2$ preserving $\mathfrak{B}$ can be lifted. Then up to equivalence, one of the following occurs:
\begin{enumerate}
  \item[\rm(1)] $A=\mathbb{Z}_{p^k}^{n-1}$, and $\phi_\pi(x_i)=\mathbf{e}_i$, $1\le i\le n-1$. 
  \item[\rm(2)] $A=\mathbb{Z}_{p^r}^{n-2}\times\mathbb{Z}_{p^k}$ for some $k>r>0$ with $p^{k-r}\mid n$, and
        $$\phi_\pi(x_i)=(\mathbf{e}_{i},1), \quad 1\le i\le n-2; \qquad \phi_\pi(x_{n-1})=(\mathbf{0},1).$$
  \item[\rm(3)] $A=\mathbb{Z}_{p^k}$ with $p^k\mid n$, and $\phi_\pi(x_i)=1$ for all $i$.
\end{enumerate}
In {\rm(1)} or {\rm(2)}, $\mathbf{e}_i$ is the vector with $1$ at the $i$-th position and $0$ elsewhere.
\end{thm}
According to \cite{AMO21} Section 3, this essentially solves the problem for abelian covers.

\bigskip

Let $G$ denote the group of automorphisms of $H_1(\Sigma_{0,n};\mathbb{Z})$ induced by all homeomorphisms of $\Sigma_{0,n}$. Clearly, $G$ is isomorphic to the permutation group on $\{x_1,\ldots,x_n\}$.
By \cite{GW17} Lemma 2.1, a homeomorphism $f$ of $\Sigma_{0,n}$ can be lifted if and only if there exists $\psi\in{\rm Aut}(A)$ such that $\psi\circ\phi_\pi=\phi_\pi\circ\alpha$, where $\alpha\in G$ is induced by $f$. This is equivalent to $\alpha(\ker\phi_\pi)=\ker\phi_\pi$, which in turn is equivalent to
\begin{align*}
\alpha(\ker\overline{\phi_\pi})=\ker\overline{\phi_\pi},
\end{align*}
where $\overline{\phi_\pi}: H_1(\Sigma_{0,n};\mathbb{Z}_{p^k})\twoheadrightarrow A$ is the map induced by $\phi_\pi$.
Instead of dealing with ${\rm Aut}(A)$, we directly work on $\ker\overline{\phi_\pi}$, reducing the problem to finding all subgroups of $H_1(\Sigma_{0,n};\mathbb{Z}_{p^k})$ that are invariant under all $\alpha\in G$.
In this viewpoint, $\pi$ is equivalent to $\pi'$ if and only if $\beta(\ker\overline{\phi_{\pi'}})=\ker\overline{\phi_\pi}$ for some $\beta\in G$.

For problems of such kind, a method was developed in \cite{CS13}. 

Some notations and conventions.
For a ring $R$, let $R^{\ell,m}$ denote the set of $\ell\times m$ matrices over $R$. 
For $X\in R^{\ell,m}$, let $X_{i,j}$ denote its $(i,j)$-entry, and let $\langle X\rangle$ denote the subgroup of $R^m$ generated by the row vectors of $X$.
For $Y\in\mathbb{Z}^{\ell,m}$, abusing the notation, we denote its image under the map $\mathbb{Z}^{\ell,m}\rightarrow\mathbb{Z}_{p^{k}}^{\ell,m}$
induced by the quotient map $\mathbb{Z}\twoheadrightarrow\mathbb{Z}_{p^{k}}$ also by $Y$.
Let $S_m$ denote the permutation group on $m$ elements, and embed it as a subgroup of ${\rm GL}(m,\mathbb{Z})$ by
identifying $\sigma\in S_m$ with the matrix (denoted by the same notation) whose $(i,j)$-entry is $\delta_{i,\sigma(j)}$, where $\delta$ is the Kronecker symbol.


\begin{lem}[\cite{CS13} Theorem 3.9]
Each subgroup $C\leqslant \mathbb{Z}_{p^k}^m$ is of the form $\langle PQ\omega\rangle$ for some $P\in\mathbb{Z}_{p^k}^{\ell,m}$, $Q\in{\rm GL}(m,\mathbb{Z})$, $\omega\in S_{m}$ such that
\begin{itemize}
  \item $0\le \ell\le m$, $P_{i,j}=\delta_{i,j}p^{r_{i}}$ with $0\le r_{1}\le\cdots\le r_{\ell}<k$;
  \item $Q_{i,i}=1$ for all $i$, and $Q_{j,i}=0\le Q_{i,j}<p^{r_{j}-r_{i}}$ for all $i<j$, where $r_{i}=k$ for $\ell<i\le m$. 
\end{itemize}
\end{lem}

Clearly, $\mathbb{Z}_{p^k}^m/C\cong\mathbb{Z}_{p^{r_\iota}}\times\cdots\times\mathbb{Z}_{p^{r_{m}}}$, with $\iota=\min\{i\colon r_i>0\}$.

\medskip

Let $b=n-1$. Take $x_1,\ldots,x_{b}$ as generators for $H_1(\Sigma_{0,n};\mathbb{Z})\cong\mathbb{Z}^{b}$. For each $\alpha\in G\cong S_{b+1}$, let $T^{\alpha}\in\textrm{GL}(b,\mathbb{Z})$ denote the matrix determined by
\begin{align*}
\alpha(x_{i})=\sum_{j=1}^{b}(T^{\alpha})_{i,j}\cdot x_{j}, \qquad i=1,\cdots,b.
\end{align*}
Obviously, if $\alpha\in S_{b}$, by which we mean $\alpha(x_{b+1})=x_{b+1}$, then $T^\alpha=\alpha$.

If $\ker\overline{\phi_\pi}=\langle PQ\omega\rangle$, then $\alpha(\ker\overline{\phi_\pi})=\langle PQ\omega T^\alpha\rangle$. Taking $\beta=\omega^{-1}\in S_b$, we have $\beta(\ker\overline{\phi_\pi})=\langle PQ\rangle$. Hence, up to equivalence, we can assume $\ker\overline{\phi_\pi}=\langle PQ\rangle$. Then $\alpha(\ker\overline{\phi_\pi})=\langle PQT^\alpha\rangle$ for each $\alpha\in G$. By the criterion given by \cite{CS13} Lemma 3.11, $\alpha(\ker\overline{\phi_\pi})=\ker\overline{\phi_\pi}$ is equivalent to
$$p^{r_j-r_i}\mid(QT^{\alpha}Q^{-1})_{i,j} \qquad \text{for\ all\ } i<j.$$

Thus, it suffices to find tuples $(r_{1},\ldots,r_b;Q)$ consisting of integers $0\le r_1\le\cdots\le r_b>0$ and a matrix $Q\in{\rm GL}(b,\mathbb{Z})$, such that
\begin{align}
Q_{i,i}=1, \quad 1\le i\le b; \qquad  Q_{j,i}=0\le Q_{i,j}<p^{r_{j}-r_{i}} \text{\ for\ all\ }i<j;   \label{eq:Q} \\
p^{r_j-r_i}\mid (QT^{\alpha}Q^{-1})_{i,j} \qquad \text{for\ all\ }i<j \text{ and\ }\alpha\in S_{b+1}.  \label{eq:criterion-1}
\end{align}
Such a tuple determines a regular branched $A$-cover $\pi$, with $A=\mathbb{Z}_{p^{r_\iota}}\times\cdots\times\mathbb{Z}_{p^{r_b}}$, $\iota=\min\{i\colon r_{i}>0\}$, and
\begin{align}
\phi_\pi(x_{i})=\big((Q^{-1})_{i,\iota},\ldots,(Q^{-1})_{i,b}\big), \qquad 1\le i\le b,  \label{eq:phi}
\end{align}
with the understanding that $\phi_\pi(x_{b+1})=-\phi_\pi(x_1)-\cdots-\phi_\pi(x_b)$.

\section{Proof of Theorem \ref{thm:main}}

As a simple observation, (\ref{eq:criterion-1}) holds if and only if
\begin{align}
p^{r_j-r_i}\mid(QXQ^{-1})_{i,j} \qquad \text{for\ all\ }i<j,   \label{eq:criterion-2}
\end{align}
for all $X$ in the subgroup of $\mathbb{Z}^{b,b}$ generated by $T^{\alpha}$, $\alpha\in S_{b+1}$.

Let $I$ denote the identity matrix. Let $E_{u}^v$ denote the matrix whose $(u,v)$-entry is $1$ and the other entries are all $0$.

For the permutation $\eta_i\in S_{b+1}$ switching $i$ and $b+1$,
$$T^{\eta_i}=\left(\begin{array}{ccccc} 1 & \ & \ & \ & \ \\ \ & \ddots & \ & \ & \ \\ -1 & \cdots & -1 & \cdots & -1 \\
 \ & \ & \ & \ddots & \ \\ \ & \ & \ & \ & 1 \end{array}\right)=I-E_i^i-\sum\limits_{j=1}^bE_i^j.$$

Clearly, for each $\alpha\in S_{b+1}\setminus S_b$, there uniquely exists $u\in\{1,\ldots,b\}$, $\sigma\in S_b$ such that $\alpha=\eta_u\sigma$.
Then
\begin{align}
T^{\sigma}-T^\alpha=(I-T^{\eta_u})T^\sigma=E_u^{\sigma^{-1}(u)}+\sum_{j=1}^{b}E_u^j.  \label{eq:T-T}
\end{align}
The difference of two such matrices can give rise to $E_{u}^{v}-E_{u}^{w}$ for any $v\ne w$.

Taking $i=1$ and $X=E_1^{v}-E_1^{w}$ in (\ref{eq:criterion-2}), we obtain
\begin{align}
(Q^{-1})_{v,j}\equiv (Q^{-1})_{w,j}\pmod{p^{r_j-r_1}}.  \label{eq:deduce}
\end{align}
In particular, setting $w=j=b-1$ and $v=b$ leads to $r_{b-1}=r_1$ (so that $r_i=r_1$ for all $i<b$). By (\ref{eq:Q}), $Q_{i,j}=0$ for all $1\le i<j<b$.

If $k=r_b=r_1$, then $Q_{i,b}=0$ for all $1\le i<b$, so that $Q=I$. In this case, $A=\mathbb{Z}_{p^k}^b$, and by (\ref{eq:phi}), $\phi(x_i)=\mathbf{e}_i$.

Now suppose $k>r_1$. Setting $w=j=b$ in (\ref{eq:deduce}) leads to
$$(Q^{-1})_{v,b}\equiv 1\pmod{p^{r_b-r_1}}.$$
Hence $Q_{v,b}\equiv -1\pmod{p^{r_b-r_1}}$ for all $v<b$.
From (\ref{eq:T-T}) we see that for $i<b$,
$$(Q(T^{\sigma}-T^\alpha)Q^{-1})_{i,b}\equiv (\delta_{u,i}-\delta_{u,b})(b+1)\pmod{p^{r_b-r_1}},$$
so $p^{r_b-r_1}\mid(Q(T^{\sigma}-T^\alpha)Q^{-1})_{i,b}$ is equivalent to
\begin{align*}
p^{r_b-r_1}\mid b+1.
\end{align*}
Once this holds, for each $\sigma\in S_b$ and all $i<b$,
$$(QT^\sigma Q^{-1})_{i,b}\equiv\sum_{v=1}^b(Q\sigma)_{i,v}\equiv 0\pmod{p^{r_b-r_1}}$$
Therefore, (\ref{eq:criterion-1}) is fulfilled.

There are two possible cases.
\begin{itemize}
  \item If $r_1>0$, then $A=\mathbb{Z}_{p^{r_1}}^{b-1}\times\mathbb{Z}_{p^k}$, and by (\ref{eq:phi}),
        $$\phi(x_i)=(\mathbf{e}_{i},1), \quad 1\le i\le b-1; \qquad \phi(x_b)=(\mathbf{0},1).$$
  \item If $r_1=0$, then $A=\mathbb{Z}_{p^k}$, and $\phi(x_i)=1$ for each $i$.
\end{itemize}

\vspace{3mm}

Haimiao Chen (orcid Id: 0000-0001-8194-1264) \\
Department of Mathematics, Beijing Technology and Business University, Beijing, China. \ \
Email: \emph{chenhm@math.pku.edu.cn}

\end{document}